\documentclass{jloganal}    

\title{On the Borel-Cantelli Lemmas, the Erd\H{o}s-R\'{e}nyi Theorem, and the Kochen-Stone Theorem}

\author{Rob Arthan}
\givenname{Rob}
\surname{Arthan}
\address{School of Electronic Engineering and Computer Science,
Queen Mary University of London}
\email{r.arthan@qmul.ac.uk}
\urladdr{}

\author{Paulo Oliva}
\givenname{Paulo}
\surname{Oliva}
\email{p.oliva@qmul.ac.uk}
\urladdr{}


\arxivreference{}
\arxivpassword{}

\volumenumber{}
\issuenumber{}
\publicationyear{}
\papernumber{}
\startpage{}
\endpage{}
\doi{}
\MR{}
\Zbl{}
\received{}
\revised{}
\accepted{}
\published{}
\publishedonline{}
\proposed{}
\seconded{}
\corresponding{}
\editor{}
\version{}

\newtheorem{thm}{Theorem}[section]    
\newtheorem{lem}[thm]{Lemma}          
\theoremstyle{definition}
\newtheorem*{rem}{Remark}             
\newcommand{\M}{\mathbb{E}}
\newcommand{\SD}{\sigma}

\def\II{\mathbb{I}}

\def\NN{\mathbb{N}}

\def\PP{\mathbb{P}}

\def\RR{\mathbb{R}}

\def\ZZ{\mathbb{Z}}

\def\Vp{\mathbf{p}}
\def\Vq{\mathbf{q}}

\def\VX{\mathbf{X}}

\begin{document}

\begin{abstract}    
In this paper we present a quantitative analysis of the first and second Borel-Cantelli Lemmas and of two of their generalisations: the Erd\H{o}s-R\'{e}nyi Theorem, and the Kochen-Stone Theorem. We will see that the first three results have \emph{direct} quantitative formulations, giving an explicit relationship between quantitative formulations of the assumptions and the conclusion. For the Kochen-Stone theorem, however, we can show that the numerical bounds of a direct quantitative formulation are not computable in general. Nonetheless, we obtain a quantitative formulation of the Kochen-Stone Theorem using Tao's notion of \emph{metastability}.
\end{abstract}

\maketitle


\section{Introduction}

Let $( A_i )_{i=1}^\infty$ be an infinite sequence of events in a probability space $(\mathcal{S}, \mathcal{E}, \PP)$. The Borel-Cantelli lemma is a classical result in probability theory, relating the convergence or divergence of the sum $\sum_{i=1}^\infty \PP[A_i]$ with the probability of the event ``$A_i$ infinitely often'', which is defined as follows:
\[ A_i~\mbox{i.o.} = \bigcap_{n=1}^\infty \bigcup_{i \geq n} A_i \]
i.e., $\omega \in \mathcal{S}$ happens infinitely often in $( A_i )_{i=1}^\infty$ if for all $n$ there exists an $i \geq n$ such that $\omega \in A_i$. Note that, as $\bigcup_{i \geq 1} A_i \supseteq \bigcup_{i \geq 2} A_i \supseteq \ldots$, we have
\[
\PP[A_i~\mbox{i.o.}] = \lim_{n \to \infty} \PP\left[ \bigcup_{i \geq n} A_i\right]
\]

The Borel-Cantelli lemma (see, for example, Feller \cite{feller68}) is normally presented in two parts. The first part says that when the sum $\sum_{i=1}^\infty \PP[A_i]$ converges, then the event $A_i~\mbox{i.o.}$ has probability zero:

\begin{thm}[First Borel-Cantelli Lemma]\label{thm:first-bcl} Let $( A_i )_{i=1}^\infty$ be an infinite sequence of events such that $\sum_{i=1}^\infty \PP[A_i] < \infty$. Then $\PP[A_i~\mbox{i.o.}] = 0$.
\end{thm}

The second part says that when $\sum_{i=1}^\infty \PP[A_i]$ diverges, and when the $A_i$ are mutually independent, then the event $A_i~\mbox{i.o.}$ has probability one:

\begin{thm}[Second Borel-Cantelli Lemma]\label{thm:second-bcl} Let $( A_i )_{i=1}^\infty$ be an infinite sequence of mutually independent events such that $\sum_{i=1}^\infty \PP[A_i] = \infty$. Then $\PP[A_i~\mbox{i.o.}] = 1$.
\end{thm}

In \cite{kochen-stone64}, Kochen and Stone presented a result that generalises the Second Borel-Cantelli Lemma in two directions: {\em(i)} it gives a lower bound
on $\PP[A_i~\mbox{i.o.}]$ when the $A_i$ are not mutually independent and {\em(ii)} it can be used to show that the assumption of mutual independence in the original lemma can be weakened to pairwise independence. We formulate this generalisation following Yan~\cite{yan06}:

\begin{thm}[Kochen-Stone] Let $( A_i )_{i=1}^\infty$ be an infinite sequence of events such that $\sum_{i=1}^\infty \PP[A_i] = \infty$. 
Then
\begin{align} 
\PP[A_i~i.o]
  & \geq \limsup_{n \to \infty} \frac{(\sum_{k=1}^n \PP[A_k])^2}{\sum_{i,k=1}^n \PP[A_i A_k]} \label{kochen-stone-ineq} 
\end{align}
\end{thm}


Erd\H{o}s and R\'{e}nyi \cite{erdos-renyi59} gave a result that is intermediate between the second Borel-Cantelli lemma and the Kochen-Stone theorem. Like the Kochen-Stone theorem it implies that the assumption of mutual independence in the second Borel-Cantelli lemma can be weakened to pairwise independence. Erd\H{o}s and R\'{e}nyi applied their theorem to the study of generalised Cantor expansions for real numbers.

\begin{thm}[Erd\H{o}s-R\'{e}nyi] \label{erdos-thm} Let $( A_i )_{i=1}^\infty$ be an infinite sequence of events such that $\sum_{i=1}^\infty \PP[A_i] = \infty$ and
\begin{align} \label{erdos-renyi-eq}
\liminf_{n \to \infty} \frac{\sum_{i,k=1}^n \PP[A_i A_k]}{(\sum_{k=1}^n \PP[A_k])^2} = 1
\end{align}
Then $\PP[A_i~i.o] = 1$.
\end{thm}

The aim of the present note is to investigate ``quantitative'' versions of each of these ``qualitative'' results. The methods we use come from the \emph{proof mining} programme (see Kohlenbach \cite{Kohlenbach(2008)}) -- where numerical information is obtained from (often non-constructive) proofs via logical methods. For some noteworthy applications of these methods, see the work of Avigad and collaborators on Ergodic Theory \cite{Avigad2009} and Kohlenbach and collaborators on Fixed-Point Theory \cite{Kohlenbach(2009), Kohlenbach(2012)}. 

Terence Tao's programme of bridging ``soft'' and ``hard'' analysis \cite{tao2008structure} was an independent rediscovery of some of these ideas. The results as presented above are results of ``soft analysis'': they relate statements about convergence or divergence, without giving any numeric information about the corresponding rates of convergence or divergence. For instance, regarding the First Borel-Cantelli Lemma, it is natural to ask how the rate of convergence of the sequence of partial sums $(\sum_{i=1}^n \PP[A_i])_{n=1}^\infty$ relates to the rate of convergence of the sequence of probabilities $(\PP[\bigcup_{i \geq n} A_i])_{n=1}^\infty$, whose limit is $\PP[A_i~\mbox{i.o.}]$. Similar questions arise regarding the other three results.

We provide here answers to these four questions. We will find in Section~\ref{sec:b-c} that the answer is almost trivial for the First Borel-Cantelli Lemma, as it has a very direct (constructive) proof. It turns out that the sequence  $(\PP[\bigcup_{i \geq n} A_i])_{n=1}^\infty$ converges with the same rate as the sequence $(\sum_{i=1}^n \PP[A_i])_{n=1}^\infty$. The answers are found to be less trivial in Section~\ref{sec:b-c} for the Second Borel-Cantelli Lemma  and in Section~\ref{sec:e-r} for the Erd\H{o}s-R\'{e}nyi theorem, but still, the quantitative versions of these follow the standard proofs of the qualitative versions quite closely. As we will see in Section~\ref{sec:k-s}, in the case of the Kochen-Stone Theorem the situation is more complicated. In Section~\ref{subsec:necessity}, we prove that a direct (computable) rate of convergence does not exist: we can find a concrete sequence of events (with computable probabilities $\PP[A_i]$) such that the rate of convergence for the quantitative version of the theorem is not computable. To allow for this, we use Tao's notion of ``rate of metastability'' (see Section~\ref{sec:tao}), a concept which is logically equivalent to convergence, but is computationally weaker. We give a quantitative version of the Kochen-Stone Theorem that provides a rate of metastability for inequality (\ref{kochen-stone-ineq}) as a computable function of the rate of divergence of the sequence of partial sums $(\sum_{i=1}^n \PP[A_i])_{n=1}^\infty$.

In Section \ref{sec-optimal} we also consider the optimality of the bounds we obtain. For the First and Second Borel-Cantelli Lemmas, we can argue that the bounds obtained are in some sense best possible. For the Erd\H{o}s-R\'{e}nyi we conjecture that the bounds we present are optimal, but do not have a proof yet. For the Kochen-Stone theorem, however, the ``metastable'' reformulation makes it much less clear what the right notion of optimality should be, and we leave this to future work.



The work presented in this paper was motivated by our on-going work in the area of metric Diophantine approximation, more specifically, on quantitative analyses of generalisations of the Khintchine-Groshev theorem on approximability of real numbers by rationals. To introduce these generalisations, let $\II^{nm}$ denote the unit cube $[0, 1]^{nm}$ in $\RR^{nm}$, and let $\psi \colon \NN \to \RR^+$ be given. A point $\VX \in \II^{nm}$, viewed as an $n \times m$ matrix, is said to be $\psi$-approximable if there are infinitely many $(\Vp, \Vq) \in \ZZ^m \times \ZZ^n$ such that $\| \Vq \VX + \Vp \| < \psi(\|\Vq\|)$ (where $\|\cdot\|$ is the supremum norm). Let ${\mathcal A}_{n, m}$ denote the set of $\psi$-approximable points $\VX \in \II^{nm}$. The generalised Khintchine-Groshev theorems are $0$-$1$ laws for the Lebesgue measure of the sets ${\mathcal A}_{n, m}$ governed by assumptions on the divergence (and possibly monotonicity) of sequences defined in terms of $\psi$. The most general theorem of this form is given in Beresnevich and Velani \cite{beresnevich-velani09}, which improves on earlier work of  Gallagher \cite{gallagher65} dealing with the case $n = 1$. In both of these works, the proofs break into two parts: A proof of a $0$-$1$ law and a proof that a certain set has positive measure, and hence measure $1$ by the $0$-$1$ law. The Borel-Cantelli lemmas and their generalisations are important tools in some of these proofs. In particular, the Kochen-Stone theorem is a key step\footnote{Beresnevich and Velani in fact refer to Sprindzuk \cite[Lemma 5]{sprindzuk79} for the result and not to Kochen and Stone \cite{kochen-stone64}. Presumably Sprindzuk obtained the result independently.} in Beresnevich and Velani's work \cite{beresnevich-velani09}. A quantitative analysis of these tools seemed to us to be a worthwhile investigation in its own right.

\subsection{Rate of convergence vs rate of metastability}\label{sec:tao}

As mentioned above, in our quantitative analysis of the Kochen-Stone theorem we will make use of Terence Tao's notion of \emph{metastability} \cite{tao2008structure}. As an example, consider the statement that a sequence of reals $(x_n)_{n=1}^\infty$ is Cauchy convergent:
\begin{equation} \label{cauchy-convergence}
\forall \ell \exists k \forall m, n > k \left(| x_m - x_n | < \frac{1}{2^\ell}\right).
\end{equation}

A \emph{rate of Cauchy convergence} (or just \emph{rate of convergence}) for the sequence is a function $\phi \colon \NN \to \NN^+$ (where $\NN^+$ denotes the positive integers) such that
\begin{equation} 
\forall \ell \forall m, n > \phi(\ell) \left(| x_m - x_n | < \frac{1}{2^\ell}\right). \label{mod-conv}
\end{equation}

While the mere existence of $\phi$ in (\ref{mod-conv}) is equivalent to (\ref{cauchy-convergence}), if one has an explicit $\phi$ for which (\ref{mod-conv}) holds, one has a quantitative rather than merely qualitative understanding of the convergence of the sequence $(x_n)_{n=1}^\infty$.
However, for a given sequence that is known to be convergent, it may not be possible to provide an explicit rate of convergence: in many interesting cases, it can be shown that no computable function $\phi$ satisfying (\ref{mod-conv}) exists. In such cases, it is often worth considering the equivalent ``metastable'' version of (\ref{cauchy-convergence}), namely:
\begin{equation} \label{cauchy-convergence-meta}
\forall \ell, f^{\NN^+\to\NN^+}\exists k \forall m, n \in [k, f(k)] \left(| x_m - x_n | < \frac{1}{2^\ell}\right).
\end{equation}
Clearly (\ref{cauchy-convergence}) directly implies (\ref{cauchy-convergence-meta}). But (\ref{cauchy-convergence-meta}) also implies (\ref{cauchy-convergence}). Assume (\ref{cauchy-convergence-meta}) and suppose (\ref{cauchy-convergence}) does not hold for some $\ell$, i.e.
\[ \forall k \exists m, n > k \left(| x_m - x_n | \geq \frac{1}{2^\ell}\right) \]
and let $f \colon \NN^+\to \NN^+$ be any function which provides an upper bound for $m$ and $n$ for each given $k$, i.e.
\[ \forall k \exists m, n \in [k, f (k)] \left(| x_m - x_n | \geq \frac{1}{2^\ell}\right) \]
Taking this $f$ in (\ref{cauchy-convergence-meta}) leads to a contradiction.

In cases where there is no computable rate of convergence (as will be the case with the Kochen-Stone theorem), one can still attempt to find a computable \emph{rate of metastability} instead. In the convergence example above, the rate of metastability would be a function $\phi \colon \NN \times (\NN^+\to \NN^+) \to \NN^+$ such that
\begin{equation} \label{cauchy-convergence-meta-rate}
\forall \ell, f^{\NN^+\to\NN^+}\exists k \leq \phi(\ell, f) \forall m, n \in [k, f (k)] \left(| x_m - x_n | < \frac{1}{2^\ell}\right).
\end{equation}
One should think of the function $f \colon \NN^+\to \NN^+$ as potentially producing longer and longer intervals $[k, f(k)]$, and the rate of metastability $\phi$ as trying to find, for arbitrarily large $\ell$, an interval in which the sequence is $\frac{1}{2^\ell}$-stable, i.e. $\forall m, n \in [k, f(k)] (| x_m - x_n | < \frac{1}{2^\ell})$.

\subsection{Rate of divergence}

Let $(x_i)_{i=1}^\infty$ be a non-decreasing sequence of real numbers. The sequence is said to \emph{diverge} if
\[ \forall N \exists i (x_i \geq N) \]
We say that a function $\omega \colon \NN^+ \to \NN^+$ gives a \emph{rate of divergence for $(x_i)_{i=1}^\infty$} if
\[ \forall N (x_{\omega(N)} \geq N) \]
In the sequel, the sequence $x_i$ will typically comprise the partial sums of a series of terms in the interval $[0, 1]$, for which we have the following lemma. 

\begin{lem} \label{divergence-lemma} Let $x_n = \sum_{i=1}^n a_i$, where $0 \leq a_i \leq 1$ for all $i$, and assume $\omega \colon \NN^+ \to \NN^+$ is a rate of divergence for $(x_n)_{n=1}^\infty$. Then, for all $n, N$ we have
\begin{equation}
    \sum_{i=n}^{\omega(n + N - 1)} a_i \geq N
\end{equation}
\end{lem}
\begin{proof} Since $\omega$ is a rate of divergence and each $a_i \le 1$, we have that
\begin{equation}
    \sum_{i=n}^{\omega(n + N - 1)} a_i = \sum_{i=1}^{\omega(n + N - 1)} a_i - \sum_{i=1}^{n - 1} a_i \geq (n + N - 1) - (n - 1) = N
\end{equation}
\end{proof}

\section{Quantitative Borel-Cantelli Lemmas}\label{sec:b-c}

Given an infinite sequence of events $( A_i )_{i=1}^\infty$, the first Borel-Cantelli lemma says that the probability of $A_i~\mbox{i.o.}$ is 0 when $\sum_{i=1}^\infty \PP[A_i]$ is finite. Here the assumption ``$\sum_{i=1}^\infty \PP[A_i]$ is finite'' is equivalent to the convergence of the sequence of partial sums $s_k = \sum_{i=1}^k \PP[A_i]$. In quantitative terms, that implies the existence of a rate of Cauchy convergence $\psi(\ell)$ for $s_k$, i.e.
\begin{equation}
	\forall \ell, m, n > \psi(\ell) \left( |s_m - s_n| < \frac{1}{2^\ell} \right)
\end{equation}
or, equivalently, the existence of a function $\phi(\ell)$  such that
\begin{equation}
	\forall \ell \forall m > \phi(\ell) \left( \sum_{i=\phi(\ell)}^m \PP[A_i] < \frac{1}{2^\ell} \right)
\end{equation}
%


\begin{thm}[First Borel-Cantelli Lemma -- Quantitative Version]\label{thm:first-bcl-q} Let $( A_i )_{i=1}^\infty$ be an infinite sequence of events. Assume that $(\sum_{i=1}^m \PP[A_i])_{m=1}^\infty$ converges with a rate of convergence $\phi \colon \NN \to \NN^+$, i.e. that for all $\ell \ge 0$ and $m > \phi(\ell)$
\[ \sum_{i=\phi(\ell)}^m \PP[A_i] \leq \frac{1}{2^\ell} \]
Then the sequence $(\PP[\bigcup_{i=1}^m A_i])_{m=1}^\infty$ converges to $0$ with the same rate, i.e. for all $\ell \ge 0$ and $m > \phi(\ell)$
\[
\PP\left[\bigcup_{i=\phi(\ell)}^m A_i\right] \leq \frac{1}{2^\ell}
\]
\end{thm}
\begin{proof} By subadditivity we have
\begin{equation} \label{first-lemma-proof}
    \PP\left[\bigcup_{i=\phi(\ell)}^m A_i\right] \leq \sum_{i=\phi(\ell)}^m \PP[A_i] \leq \frac{1}{2^\ell} 
\end{equation}
for all $\ell > 0$ and $m > \phi(\ell)$. \end{proof}

The second Borel-Cantelli lemma says that, under the extra assumption that the events are mutually independent, the probability of $A_i~\mbox{i.o.}$ is 1 whenever the sum $\sum_{i=1}^\infty \PP[A_i]$ diverges. In our quantitative version of this lemma we will estimate, for each $n$, how fast the sequence
\[
\left( \PP\left[\bigcup_{i=n}^m A_i\right] \right)_{m=1}^\infty
\]
\emph{converges} to 1, given a rate of \emph{divergence} for the sequence $( \sum_{i=n}^m \PP[A_i])_{m=1}^\infty$.

\begin{thm}[Second Borel-Cantelli Lemma -- Quantitative Version]\label{thm:second-bcl-q} Let $( A_i )_{i=1}^\infty$ be an infinite sequence of events which are mutually independent. Assume that the sequence $( \sum_{i=1}^n \PP[A_i])_{n=1}^\infty$ diverges with rate $\omega \colon \NN^+\to \NN^+$, i.e. for all $N$
\[ \sum_{i=1}^{\omega(N)} \PP[A_i] \geq N \]
then, for all $n$ and $N$,
\[
\PP\left[\bigcup_{i=n}^{\omega(n + N - 1)} A_i\right] \geq 1 - e^{-N}
\]
\end{thm}
\begin{proof} Fix $n$ and $N$. Let us write $\overline{A}_i$ for the complement of the event $A_i$. The independence of the events implies
\begin{align}
\PP\left[\bigcap_{i=n}^{\omega(n + N - 1)} \overline{A}_i \right] 
    & = \prod_{i=n}^{\omega(n + N - 1)} \PP \left[ \, \overline{A}_i \right] \\
    & = \prod_{i=n}^{\omega(n + N - 1)} (1 - \PP \left[ A_i \right])
\end{align}
Taking the natural logarithm on both sides we have
\begin{align}   
\ln \left( \PP\left[\bigcap_{i=n}^{\omega(n + N - 1)} \overline{A}_i \right] \right) 
    & = \ln\left(\prod_{i=n}^{\omega(n + N - 1)} (1 - \PP \left[ A_i \right])\right) \\
    & = \sum_{i=n}^{\omega(n + N - 1)} \ln(1 - \PP \left[ A_i \right]) \\
    & \leq - \sum_{i=n}^{\omega(n + N - 1)} \PP[A_i] \label{ln-lemma-use} \\
    & \leq - N \label{div-lemma-use}
\end{align}
where inequality (\ref{ln-lemma-use}) follows from the fact that $\ln(1 + x) \leq x$, for all $x \in (-1, \infty)$ and inequality (\ref{div-lemma-use}) follows from Lemma \ref{divergence-lemma}. Hence,
\begin{equation} \label{2nd-bc-proof-eq1}
\PP\left[\bigcap_{i=n}^{\omega(n + N - 1)} \overline{A}_i \right] \leq e^{-N}
\end{equation}
and so
\begin{equation} \label{2nd-bc-proof-eq2}
\PP\left[\bigcup_{i=n}^{\omega(n + N - 1)} A_i\right] \geq 1 - e^{-N}
\end{equation}
\end{proof}

\subsection{Proving qualitative version from quantitative one}

That Theorem~\ref{thm:first-bcl} (First Borel-Cantelli lemma) follows from Theorem~\ref{thm:first-bcl-q} (its quantitative version) is clear. Let us show that Theorem~\ref{thm:second-bcl} (Second Borel-Cantelli lemma) also follows from Theorem~\ref{thm:second-bcl-q} (its quantitative version).

Let $( A_i )_{i=1}^\infty$ be an infinite sequence of mutually independent events. Assuming  $\sum_{i=1}^\infty \PP[A_i] = \infty$, there exists $\omega \colon \NN^+\to \NN^+$, such that for all $N$,
\[ \left(\sum_{i=1}^{\omega(N)} \PP[A_i]\right) \geq N \]
By Theorem~\ref{thm:second-bcl-q}, for all $n$ and $N$
\[
\PP\left[\bigcup_{i=n}^{\omega(n + N - 1)} A_i\right] \geq 1 - e^{-N}
\]
which implies that
\[
\PP\left[\bigcup_{i=n}^{\infty} A_i\right] \geq 1 - e^{-N}
\]
and hence $\PP[A_i~\mbox{i.o.}] = 1$.

\section{Quantitative Erd\H{o}s-R\'{e}nyi Theorem} \label{sec:e-r}

In this section, we present a quantitative version of the Erd\H{o}s-R\'{e}nyi theorem. Our proof follows that of Erd\H{o}s and R\'{e}nyi \cite{erdos-renyi59}, but uses more modern notation: $X, Y$ for random variables, $\M(X)$ for the expectation of $X$ (or mean value in Erd\H{o}s and R\'{e}nyi's terminology), and $\SD(X)$ for the standard deviation of $X$. 

A simplistic logical formalisation of equation (\ref{erdos-renyi-eq}), involving $\liminf$, would give a formula with quantifier prefix $\forall\exists\forall\exists$. We can simplify this to a $\forall \exists$ using the following lemmas:

\begin{lem} \label{erdos-renyi-lemma} For any sequence of events $( A_i )_{i=1}^\infty$ we have that for all $n \geq 1$
\begin{equation}
\frac{\sum_{i,k=1}^n \PP[A_i A_k]}{ (\sum_{k=1}^n \PP[A_k])^2} \geq 1
\end{equation}
\end{lem}
\begin{proof} Let $X_i$ be the random variable given by the indicator function of the event $A_i$. Then $\M(X_i) = \PP[A_i]$ and $\M(X_i X_k) = \PP[A_i A_k]$. Define $Y_n = \sum_{i=1}^n X_i$ so that
\begin{equation} \label{erdos-step3}
    \frac{\sum_{i,k=1}^n \PP[A_i A_k]}{(\sum_{i=1}^n \PP[A_i])^2}
    = 
    \frac{\M(Y_n^2)}{\M^2(Y_n)}
\end{equation}
Since the variance $\SD^2(X) = \M(X^2) - \M^2(X) \ge 0$ for any $X$, the result follows.
\end{proof}

\begin{lem}\label{lma:liminf} For any sequence of reals $d_n \geq 1$, the following are equivalent:
\begin{align}
    & \liminf_{n \to \infty} d_n = 1 \\
    & \forall \ell, n \exists k \geq n \left( d_k \leq 1 + 2^{-\ell} \right)
\end{align}
\end{lem}
\begin{proof} Because $d_n \geq 1$, we have that $\liminf_{n \to \infty} d_n = 1$ is equivalent to
\[
\forall \ell \exists m \forall i \geq m \left( \inf_{k \geq i} d_k \leq 1 + 2^{-\ell} \right) 
\]
which, by the definition of $\inf$, is equivalent to
\[
\forall \ell \exists m \forall i \geq m \exists k \geq i (d_k \leq 1 + 2^{-\ell}) 
\]
This is easily seen to be equivalent to $\forall \ell, n \exists k \geq n \left( d_k \leq 1 + 2^{-\ell} \right)$ (for the right-to-left direction take $i = \max(m, n)$). \end{proof}

In the quantitative version of the Erd\H{o}s-R\'{e}nyi Theorem, we will assume that we are given a rate of divergence for the sequence $( \sum_{i=n}^m \PP[A_i])_{m=1}^\infty$, and a function $\phi$ witnessing 
\[ 
\forall \ell, n \exists m \geq n \left(
\frac{\sum_{i,k=1}^m \PP[A_i A_k]}{(\sum_{k=1}^m \PP[A_k])^2} \leq 1 + \frac{1}{2^\ell}
\right)
\]
which, by Lemmas \ref{erdos-renyi-lemma} and \ref{lma:liminf}, is equivalent to the assumption %
\[ 
\liminf_{n \to \infty} \frac{\sum_{i,k=1}^n \PP[A_i A_k]}{(\sum_{k=1}^n \PP[A_k])^2} = 1
\]

\begin{thm}[Erd\H{o}s-R\'{e}nyi --- Quantitative Version] \label{erdos-renyi-thm} Let $( A_i )_{i=1}^\infty$ be an infinite sequence of events. Let $\omega \colon \NN^+ \to \NN^+$ be such that for all $N$
\begin{equation} \label{erdos-assumption1}
    \left(\sum_{i=1}^{\omega(N)} \PP[A_i]\right) \geq N
\end{equation} 
and let $\phi \colon \NN \times \NN \to \NN$ be such that
\begin{equation} \label{erdos-assumption2}
    \forall \ell, n \left( \phi(\ell, n) \geq n \wedge \frac{\sum_{i,k=1}^{\phi(\ell, n)} \PP[A_i A_k]}{(\sum_{i=1}^{\phi(\ell, n)} \PP[A_i])^2} \leq 1 + \frac{1}{2^\ell} \right)
\end{equation}
Define $n_1 = \phi(1, 1)$ and, for $k > 1$, $n_k = \phi(k, \max(n_{k-1}, k))$. Then, for all $n$ and $\ell$
\begin{equation} \label{erdos-conclusion}
    \PP\left[\bigcup_{i=n}^{n_m} A_i\right] \geq 1 - \frac{1}{2^\ell} 
\end{equation}
where $m = \max(\omega(2n), \ell+3)$.
\end{thm}
\begin{proof} Let $X_i$ and $Y_n$ be as in the proof of Lemma \ref{erdos-renyi-lemma}. Assumption (\ref{erdos-assumption2}) gives us
\begin{equation} \label{erdos-step4}
    \forall \ell, n \left( \phi(\ell, n) \geq n \wedge 
    \frac{\M(Y_{\phi(\ell, n)}^2)}{\M^2(Y_{\phi(\ell, n)})}
    \leq 1 + \frac{1}{2^\ell} \right)
\end{equation}
Since $\M(Y_n^2) = \SD^2(Y_n) + \M^2(Y_n)$, we have
\begin{equation} \label{erdos-step5}
    \frac{\M(Y_n^2)}{\M^2(Y_n)}
    = 
    \frac{\SD^2(Y_n)}{\M^2(Y_n)}
    + 1
\end{equation}
By (\ref{erdos-step4}) and (\ref{erdos-step5}), we get
\begin{equation} \label{erdos-step6}
    \forall \ell, n \left( \phi(\ell, n) \geq n \wedge 
    \frac{\SD^2(Y_{\phi(\ell, n)})}{\M^2(Y_{\phi(\ell, n)})}
    \leq \frac{1}{2^\ell} \right)
\end{equation}
Let $n_1 = \phi(1, 1)$ and, for $k > 1$, $n_k = \phi(k, \max(n_{k-1}, k))$. The above implies (taking $\ell =k$ and $n = \max(n_{k-1}, k)$)
\begin{equation} \label{erdos-step7}
    \forall k \left( 
    \frac{\SD^2(Y_{n_k})}{\M^2(Y_{n_k})}
    \leq \frac{1}{2^k} \right)
\end{equation}
The Chebyshev inequality tells us that
\begin{equation} \label{chebyshev}
    \PP[|Y_n - \M(Y_n)| \geq \lambda \SD(Y_n)] \leq \frac{1}{\lambda^2}
\end{equation}
Taking $\lambda = \frac{\varepsilon \M(Y_n)}{\SD(Y_n)}$, we have that, for any given $\varepsilon \in (0,1)$,
\begin{equation} \label{erdos-step8}
    \PP[Y_n \leq (1-\varepsilon) \M(Y_n)]
    \leq
    \frac{\SD^2(Y_n)}{\varepsilon^2 \M^2(Y_n)}
\end{equation}
From (\ref{erdos-step7}) and (\ref{erdos-step8}) (taking $n = n_k$ and $\varepsilon = 1/2$), we find that
\begin{equation} \label{erdos-step9}
    \forall k \left( \PP\left[Y_{n_k} \leq \frac{\M(Y_{n_k})}{2}\right] 
    \leq
    \frac{1}{2^{k-2}} \right)
\end{equation}
Let $B_k$ be the event $Y_{n_k} \leq \frac{\M(Y_{n_k})}{2}$, so that (\ref{erdos-step9}) together with the formula for the partial sums of a geometric series implies
\begin{equation} \label{erdos-step11}
    \forall \ell \forall m \geq \ell + 3 \left( \sum_{k=\ell+3}^m \PP[B_k] 
    \leq
    \frac{1}{2^\ell} \right)
\end{equation}
By the quantitative version of the First Borel-Cantelli Lemma (Theorem \ref{thm:first-bcl-q}) 
\begin{equation} \label{erdos-step12}
    \forall \ell \forall m \geq \ell + 3 \left( \PP\left[\bigcup_{k=\ell+3}^m B_k\right] 
    \leq
    \frac{1}{2^\ell} \right)
\end{equation}
Hence
\begin{equation} \label{erdos-step13}
    \forall \ell \forall m \geq \ell + 3 \left( \PP\left[\bigcap_{k=\ell+3}^m \overline{B}_k \right] 
    \geq
    1 - \frac{1}{2^\ell} \right)
\end{equation}
where $\overline{B}_k$ is the complement of $B_k$, i.e. the event $Y_{n_k} > \frac{\M(Y_{n_k})}{2}$. Hence, an outcome $x \in \overline{B}_k$ if $Y_{n_k}(x) = \sum_{i=1}^{n_k} X_i(x) > \frac{\M(Y_{n_k})}{2}$. Therefore, $x \in \overline{B}_k$ iff $x \in A_i$ for at least $\frac{\M(Y_{n_k})}{2}$ values of $i \in [1, n_k]$. \\[1mm]
Fix $n$ and $\ell$. We will show that
\begin{equation}
    \PP\left[\bigcup_{i=n}^{n_m} A_i\right] \geq 1 - \frac{1}{2^\ell} \label{erdos-step16}
\end{equation}
for $m = \max(\omega(2n), \ell+3)$. Let $C_\ell^m = \bigcap_{k=\ell+3}^m \overline{B}_k$. By (\ref{erdos-step13}), inequality (\ref{erdos-step16}) will follow if we can show that 
\begin{equation}
    C_\ell^m \subseteq \bigcup_{i=n}^{n_m} A_i \label{erdos-step15}
\end{equation}
So, let us assume that $x \in C_\ell^m$. Since, $x \in C_\ell^m$ iff $\forall k \in [\ell + 3, m] (x \in \overline{B}_k)$, we have, taking $k = m$,
\begin{equation}
Y_{n_m}(x) > \frac{\M(Y_{n_m})}{2} 
\end{equation}
By the definition of $n_k$ and the assumption (\ref{erdos-assumption2}), which implies $\phi(\ell, n) \geq n$, we have $n_{\omega(2n)} \geq \omega(2n)$, and by the definition of $m$ we also have $m \geq \omega(2n)$. Hence
\begin{equation} \label{erdos-step14}
Y_{n_m}(x) >
\frac{\M(Y_{n_m})}{2} \geq  \frac{\M(Y_{\omega(2n)})}{2} \geq n
\end{equation}
using, for the last inequality, that assumption (\ref{erdos-assumption1}) implies $\forall N \left( \M(Y_{\omega(N)}) \geq N \right)$. 
Since, $Y_{n_m}(x) = \sum_{i=1}^{n_m} X_i(x)$, (\ref{erdos-step14}) and the pigeon-hole principle imply that for at least one $i \in [n, n_m]$ we have $X_i(x) = 1$, i.e. $x \in A_i \subseteq \bigcup_{j=n}^{n_m} A_j$.
Since $x$ was an arbitrary element of $C_\ell^m$, this gives us (\ref{erdos-step15}) and hence (\ref{erdos-step16}).
\end{proof}

\begin{rem} \label{erdos-renyi-remark}
To obtain a quantitative version of the Erd\H{o}s-R\'{e}nyi theorem, we have had to make two choices about points that are left open in the qualitative proof of \cite{erdos-renyi59}. The first choice is essentially forced upon us by our decision to use $\frac{1}{2^\ell}$ in the formulation of assumption (\ref{erdos-assumption2}). This means that where Erd\H{o}s and R\'{e}nyi take the $n_k$ to be any sequence such that 
$\sum_{k = 1}^\infty \frac{\SD^2(Y_{n_k})}{\M^2(Y_{n_k})}$ converges, we have had to choose a sequence such that the series is dominated by a geometric series (see formula (\ref{erdos-step7})).
The second choice is that in formula (\ref{erdos-step9}), we have taken $\varepsilon$ to be $\frac{1}{2}$, where Erd\H{o}s and R\'{e}nyi leave it  unspecified. As pointed out by one of the referees, both the statements and the proofs above (Lemma \ref{lma:liminf} and Theorem \ref{erdos-renyi-thm}) could be made more general, but rather more complicated, by introducing a convergent series in place of $\frac{1}{2^l}$ and a constant $\varepsilon \in (0, 1)$ in place of $\frac{1}{2}$ as parameters.
\end{rem}

\subsection{Proving qualitative version from quantitative one}

Let us show that the Erd\H{o}s-R\'{e}nyi theorem (Theorem \ref{erdos-thm}) follows directly from our quantitative version (Theorem \ref{erdos-renyi-thm}).

Let $( A_i )_{i=1}^\infty$ be an infinite sequence of events such that $\sum_{i=1}^\infty \PP[A_i] = \infty$ and assume
\begin{align} 
\liminf_{n \to \infty} \frac{\sum_{i,k=1}^n \PP[A_i A_k]}{(\sum_{k=1}^n \PP[A_k])^2} = 1
\end{align}
The above assumptions imply that there exists an $\omega \colon \NN^+ \to \NN^+$ such that
\begin{equation}
    \forall N \left( \sum_{i=1}^{\omega(N)} \PP[A_i] \geq N \right)
\end{equation} 
and a function $\phi \colon \NN \times \NN \to \NN$ such that
\begin{equation}
    \forall \ell, n \left( \phi(\ell, n) \geq n \wedge \frac{\sum_{i,k=1}^{\phi(\ell, n)} \PP[A_i A_k]}{(\sum_{i=1}^{\phi(\ell, n)} \PP[A_i])^2} \leq 1 + 2^{-\ell} \right)
\end{equation}
From these and Theorem \ref{erdos-renyi-thm} we have (\ref{erdos-conclusion}), which implies $\PP[A_i~i.o] = 1$.

\section{Quantitative Kochen-Stone Theorem}\label{sec:k-s}


As with the quantitative version of the second Borel-Cantelli lemma, we will also assume that we are given a rate of divergence for the sequence $( \sum_{i=n}^m \PP[A_i])_{m \in \NN^+}$. Our quantitative version will follow closely the very concise proof of the Kochen-Stone theorem discovered by Yan \cite{yan06}.

After expressing $\PP[A_i~\mbox{i.o.}]$ as a limit, the Kochen-Stone inequality (\ref{kochen-stone-ineq}) has the form $\lim_{n \to \infty}p_n \geq \limsup_{n \to \infty} q_n$. Much as in Lemma \ref{lma:liminf}, we can be more economical with the quantifiers using the following lemma:

\begin{lem}\label{lma:io-ge-limsup} For any sequence of events $( A_i )_{i=1}^\infty$ and sequence of reals $( x_i )_{i=1}^\infty$, the following are equivalent:
\begin{align}
    & \PP[A_i~i.o]\geq \limsup_{i \to \infty} x_i \label{equiv-lm-1} \\
    & \forall m, \ell \exists n > m \forall j > n \left(\PP\left[\bigcup_{i=m+1}^n A_i\right] + \frac{1}{2^\ell} \geq x_j\right) \label{equiv-lm-2}
\end{align}
\end{lem}
\begin{proof} By the definition of $A_i~\mbox{i.o.}$, (\ref{equiv-lm-1}) is equivalent to
\begin{equation} \label{equiv-lm-pf-1}
    \forall m \left( \PP\left[\bigcup_{i=m+1}^\infty A_i\right] \geq \limsup_{i \to \infty} x_i \right)    
\end{equation}
Let us first show that the above is equivalent to
\begin{equation} \label{equiv-lm-pf-2}
    \forall m , \ell \exists n > m \left( \PP\left[\bigcup_{i=m+1}^n A_i\right] + \frac{1}{2^\ell} \geq \sup_{i>n} x_i \right)    
\end{equation}
Assume (\ref{equiv-lm-pf-1}), and fix $m$ and $\ell$. Pick $n_1 > m$ such that
\begin{equation} \label{equiv-lm-pf-3}
\PP\left[\bigcup_{i=m+1}^{n_1} A_i\right] + \frac{1}{2^{\ell+1}}
>
\PP\left[\bigcup_{i=m+1}^\infty A_i\right]
\end{equation}
and $n_2 > m$ such that
\begin{equation} \label{equiv-lm-pf-4}
\limsup_{i \to \infty} x_i + \frac{1}{2^{\ell+1}}
>
\sup_{i > n_2} x_i
\end{equation}
Then, taking $n = \max(n_1, n_2)$, by (\ref{equiv-lm-pf-1}), (\ref{equiv-lm-pf-3}) and (\ref{equiv-lm-pf-4}) we get
\begin{align}
    \PP\left[\bigcup_{i=m+1}^n A_i\right] + \frac{1}{2^\ell}
        & > \PP\left[\bigcup_{i=m+1}^\infty A_i\right] + \frac{1}{2^{\ell+1}}  \geq \limsup_{i \to \infty} x_i + \frac{1}{2^{\ell+1}}  > \sup_{i>n} x_i 
\end{align}
Thus (\ref{equiv-lm-pf-1}) implies (\ref{equiv-lm-pf-2}).
Now, suppose (\ref{equiv-lm-pf-2}) holds but (\ref{equiv-lm-pf-1}) does not. Then, for some $m$ and $\ell$, we have that
\begin{equation}\label{equiv-lm-pf-5}
    \PP\left[\bigcup_{i=m+1}^\infty A_i\right] + \frac{1}{2^\ell} < \limsup_{i \to \infty} x_i
\end{equation}
But by (\ref{equiv-lm-pf-2}) we have an $n > m$ such that 
\begin{equation}
    \PP\left[\bigcup_{i=m+1}^n A_i\right] + \frac{1}{2^\ell} \geq \sup_{i>n} x_i
\end{equation}
and hence
\begin{equation}
    \PP\left[\bigcup_{i=m+1}^\infty A_i\right] + \frac{1}{2^\ell} \geq \PP\left[\bigcup_{i=m+1}^n A_i\right] + \frac{1}{2^\ell} \geq \sup_{i>n} x_i \geq \limsup_{i \to \infty} x_i
\end{equation}
contradicting (\ref{equiv-lm-pf-5}). Thus (\ref{equiv-lm-pf-2}) implies, and hence is equivalent to (\ref{equiv-lm-pf-1}). Finally, (\ref{equiv-lm-pf-2}) is equivalent to (\ref{equiv-lm-2}), by the definition of $\sup$. \end{proof}

Using the above lemma we can show that the Kochen-Stone inequality (\ref{kochen-stone-ineq}) can be equivalently written as a $\forall \exists \forall$-statement:

\begin{lem} The Kochen-Stone inequality (\ref{kochen-stone-ineq}) is equivalent to
\begin{equation} \label{kochen-stone-eq-4}
\forall m, \ell \exists n > m \forall j > n \left( \PP\left[\bigcup_{i=m+1}^n A_i\right] + \frac{1}{2^\ell} \geq \frac{(\sum_{k=1}^j \PP[A_k])^2}{\sum_{i,k=1}^j \PP[A_i A_k]} \right)
\end{equation}
\end{lem}
\begin{proof} This follows from Lemma~\ref{lma:io-ge-limsup} using the sequence $x_j = \frac{(\sum_{k=1}^j \PP[A_k])^2}{\sum_{i,k=1}^j \PP[A_i A_k]}$.
\end{proof}

As we will see in Section \ref{meta-necessity}, there can be no computable function of $m$ and $\ell$ that bounds $n$ in (\ref{kochen-stone-eq-4}). Therefore, we will consider its metastable counterpart:
\begin{equation} \label{kochen-stone-eq-5}
\forall m, \ell, g \exists n > m \forall j \in [n, g (n)] \left( \PP\left[\bigcup_{i=m+1}^n A_i\right] + \frac{1}{2^\ell} \geq \frac{(\sum_{k=1}^j \PP[A_k])^2}{\sum_{i,k=1}^j \PP[A_i A_k]} \right)
\end{equation}
and will produce an explicit computable bound on $n$ as a function of $m, \ell$ and the function $g$. 

\begin{thm}[Kochen-Stone --- Quantitative Version] \label{quant-kochen-stone-thm} Let $( A_i )_{i=1}^\infty$ be an infinite sequence of events. Let $\omega \colon \NN^+ \to \NN^+$ be such that for all $N$
\[ \left(\sum_{i=1}^{\omega(N)} \PP[A_i]\right) \geq N \]
Then, for all $m$ and $\ell$ and $g \colon \NN^+ \to \NN^+$ such that $g(i) > i$, for all $i$, there exists an $n > m$ such that
\begin{itemize}
    \item $n \leq g^{(2^{\ell+1})}(\max(\omega(2^{\ell+2} \sum_{i=1}^m \PP[A_i]), m))$, and
    \item for all $j \in [n , g(n)]$
\begin{equation}
\PP\left[\bigcup_{i=m+1}^n A_i\right] + \frac{1}{2^\ell} \geq \frac{(\sum_{i=1}^j \PP[A_i])^2}{\sum_{i,k=1}^j \PP[A_i A_k]} 
\end{equation}
\end{itemize}
%
\end{thm}

\begin{rem} Since $\PP[A_i] \leq 1$, we have that $\sum_{i=1}^m \PP[A_i] \leq m$. Hence, we can obtain a bound $g^{(2^{\ell+1})}(\omega(2^{\ell + 2} m))$ on $n$ which is completely independent of the actual events $A_i$, but only depends on the parameters $\omega, g, m$ and $\ell$.
\end{rem}

Before we embark on the proof of Theorem \ref{quant-kochen-stone-thm}, we need three further lemmas.

\begin{lem}[Chung-Erd\H{o}s Inequality] \label{lemma-1} For every $n$ and sequence of events $A_1, \ldots, A_n$
\[
\PP \left[ \bigcup_{k=1}^n A_k \right] \geq \frac{(\sum_{k=1}^n \PP[A_k])^2}{\sum_{i,k=1}^n \PP[A_i A_k]}
\]
\end{lem}
\begin{proof} See \cite{chung-erdos52} or \cite{yan06}. \end{proof}

\begin{lem} \label{kochen-stone-lemma} Let $a$ and $b$ be such that $0 < a \leq b$. Assume $0 \leq x < a, 0 \leq y < b$ and $0 < \varepsilon$. If $b \geq 4 x / \varepsilon^2$ then
\[ 
\frac{(\sqrt{a} - \sqrt{x})^2}
{b - y} + \varepsilon \geq \frac{a}{b} 
\]
\end{lem}
\begin{proof} Since
\[\frac{(\sqrt{a} - \sqrt{x})^2}
{b - y} = \frac{a - 2 \sqrt{x a} + x}{b - y} \geq \frac{a- 2 \sqrt{x a}}{b}
\]
it is enough to show that
\[
\frac{2\sqrt{x a}}{b} \leq \varepsilon
\]
But since $a \leq b$, that follows from
\[
\frac{2\sqrt{x b}}{b} = 2\frac{\sqrt{x}}{\sqrt{b}} \leq \varepsilon
\]
which follows from $b \geq 4 x / \varepsilon^2$. \end{proof}

\begin{lem}\label{lma:k-s-estimate}
Let $( A_i )_{i=1}^\infty$ and $\omega \colon \NN^+ \to \NN^+$ be as in the statement of Theorem~\ref{quant-kochen-stone-thm}.
Then for all $m$ and $\varepsilon > 0$ and all $j > \max\left(\omega\left(\left\lceil\frac{2 \sum_{i=1}^m \PP[A_i]}{\varepsilon}\right\rceil\right), m\right)$
\[
\PP\left[\bigcup_{i=m+1}^j A_i\right] + \varepsilon \geq \frac{(\sum_{i=1}^j \PP[A_i])^2}{\sum_{i,k=1}^j \PP[A_i A_k]} 
\]
\end{lem}
\begin{proof} Let $m$ and $\varepsilon > 0$ be fixed. Let\footnote{In the notation of the proof of Lemma \ref{erdos-renyi-lemma}, $a_n = \M^2(Y_n)$ and $b_n = \M(Y_n^2)$.} $a_n = (\sum_{i=1}^n \PP[A_i])^2$ and $b_n = \sum_{i,k=1}^n \PP[A_iA_k]$. The assumption in the statement of Theorem~\ref{quant-kochen-stone-thm} says that $\sqrt{a_n}$ diverges with rate $\omega$. By Lemma \ref{lemma-1}, $b_n \geq a_n$. Hence, for all $N$,
\[
b_{\omega(N)} \geq a_{\omega(N)} \geq N^2 
\]
Since $\sum_{i,k=m+1}^n \PP[A_i A_k] \leq b_n - b_m$, by Lemma \ref{lemma-1}, we have for all $j > m$,
\begin{equation} \label{lma:k-s-estimate-step1}
\PP \left[ \bigcup_{i=m+1}^j A_i \right] + \varepsilon
\geq 
\frac{(\sqrt{a_j} - \sqrt{a_m})^2}{b_j - b_m} + \varepsilon 
\end{equation}
Let $M = \max(\omega\left(\left\lceil\frac{2 \sum_{i=1}^m \PP[A_i]}{\varepsilon}\right\rceil\right), m)$. And let $j > M$. By assumption we have that
\[ b_j \geq b_{\omega\left(\left\lceil\frac{2 \sum_{i=1}^m \PP[A_i]}{\varepsilon}\right\rceil\right)} \geq \left(\frac{2 \sum_{i=1}^m \PP[A_i]}{\varepsilon}\right)^2 \geq \frac{4 a_m}{\varepsilon^2} \]
By Lemma \ref{kochen-stone-lemma}, we have that 
\begin{equation} \label{lma:k-s-estimate-step2}
\frac{(\sqrt{a_j} - \sqrt{a_m})^2}{b_j - b_m} + \varepsilon 
\geq 
\frac{a_j}{b_j} = \frac{(\sum_{i=1}^j \PP[A_i])^2}{\sum_{i,k=1}^j \PP[A_i A_k]}
\end{equation}
The result follows from (\ref{lma:k-s-estimate-step1}) and (\ref{lma:k-s-estimate-step2}).
\end{proof}

\begin{proof}[Proof of Theorem \ref{quant-kochen-stone-thm}] Let $n_0 = \max\left(\omega(\lceil2^{\ell+2} \sum_{i=1}^m \PP[A_i]\rceil), m\right)$ and $n_{r+1} = g(n_r)$. We claim that with $n = n_r$ for some $r \leq 2^{\ell+1}$, the conclusion holds. Assume this is not the case. Then for each $r$ there is a $j_r \in [n_r, n_{r+1}]$ such that 
\begin{equation}
    \PP\left[\bigcup_{i=m+1}^{n_r} A_i\right] + \frac{1}{2^\ell} < 
    \frac{(\sum_{i=1}^{j_r} \PP[A_i])^2}{\sum_{i,k=1}^{j_r} \PP[A_i A_k]}
    \label{ks-proof-form-1}
\end{equation}
But by Lemma~\ref{lma:k-s-estimate}, with $\varepsilon = \frac{1}{2^{\ell+1}}$
\begin{equation}
    \frac{(\sum_{i=1}^{j_r} \PP[A_i])^2}{\sum_{i,k=1}^{j_r} \PP[A_i A_k]}
\leq
\PP\left[\bigcup_{i=m+1}^{j_r} A_i\right] + \frac{1}{2^{\ell+1}}
    \label{ks-proof-form-2}
\end{equation}
and hence, combining (\ref{ks-proof-form-1}) and (\ref{ks-proof-form-2}) and subtracting $\frac{1}{2^{\ell+1}}$, we have:
\begin{equation}
    \PP\left[\bigcup_{i=m+1}^{n_r} A_i\right] + \frac{1}{2^{\ell+1}} <
    \PP\left[\bigcup_{i=m+1}^{j_r} A_i\right]
\le
\PP\left[\bigcup_{i=m+1}^{n_{r+1}} A_i\right]
\label{ks-proof-form-3}
\end{equation}
%
Chaining together the inequalities given by (\ref{ks-proof-form-3}) for $r= 0$ to $2^{\ell+1}$, we have
\begin{equation}
    \PP\left[\bigcup_{i=m+1}^{n_0} A_i\right] + 1 
<
\PP\left[\bigcup_{i=m+1}^{n_{2^{\ell+1}}} A_i\right]
\leq 1
\end{equation}
which is a contradiction.
\end{proof}

\subsection{Proving qualitative version from quantitative one}

Let us argue that the qualitative version of the Kochen-Stone theorem (Theorem \ref{quant-kochen-stone-thm}) directly implies the original qualitative version. Theorem \ref{quant-kochen-stone-thm} implies that for all $m, \ell$ and $g \colon \NN^+ \to \NN^+$ (with $g(i) > i$) there exists an $n$ such that
\[
\forall j \in [n, g(n)] \left(
\PP\left[\bigcup_{i=m+1}^{n} A_i\right] + \frac{1}{2^\ell} \geq \frac{(\sum_{k=1}^j \PP[A_k])^2}{\sum_{i,k=1}^j \PP[A_i A_k]} 
\right)
\]
But (see Section \ref{sec:tao}), the above is equivalent to: for all $m, \ell$ there exists an $n$ such that
\[
\forall j \geq n \left(
\PP\left[\bigcup_{i=m+1}^{n} A_i\right] + \frac{1}{2^\ell} \geq \frac{(\sum_{k=1}^j \PP[A_k])^2}{\sum_{i,k=1}^j \PP[A_i A_k]} 
\right)
\]
Then, this is equivalent to: For all $m, \ell$ there exists an $n$ such that
\[
\PP\left[\bigcup_{i=m+1}^{n} A_i\right] + \frac{1}{2^\ell} \geq \sup_{j \geq n} \frac{(\sum_{k=1}^j \PP[A_k])^2}{\sum_{i,k=1}^j \PP[A_i A_k]} 
\]
Hence, for all $m$ and $\ell$ there exists an $n$ such that
\[
\PP\left[\bigcup_{i=m+1}^{n} A_i\right] + \frac{1}{2^\ell}
    \geq \sup_{j \geq n} \frac{(\sum_{k=1}^j \PP[A_k])^2}{\sum_{i,k=1}^j \PP[A_i A_k]} 
    \geq \limsup_{n \to \infty} \left( \frac{(\sum_{k=1}^n \PP[A_k])^2}{\sum_{i,k=1}^n \PP[A_i A_k]} \right)
\]
Therefore, for all $m$ 
\[
\PP\left[\bigcup_{i=m+1}^{\infty} A_i\right] 
    \geq \limsup_{n \to \infty} \left( \frac{(\sum_{k=1}^n \PP[A_k])^2}{\sum_{i,k=1}^n \PP[A_i A_k]} \right)
\]
and hence 
\[ 
\PP[A_i~i.o] \geq \limsup_{n \to \infty} \left( \frac{(\sum_{k=1}^n \PP[A_k])^2}{\sum_{i,k=1}^n \PP[A_i A_k]} \right)
\]

\subsection{Necessity for use of metastability}\label{subsec:necessity}
\label{meta-necessity}

We wish to show that there is no effective bound on the witness $n$ in (\ref{kochen-stone-eq-4}), so that the approach via metastability is necessary. To this end, we will need examples where the Kochen-Stone inequality (\ref{kochen-stone-ineq}) is actually an equality with $\PP[A_i~\mbox{i.o.}] < 1$. To do this we will use the following result of Yan which shows that the diagonal terms in the sums on the right-hand side of the inequality are negligible. Yan's sketch of the proof
in \cite{yan06} is very terse, so we give more detail here.
 
\begin{thm} \label{thm:yan} Let $( A_i )_{i=1}^\infty$ be an infinite sequence of events such that $\sum_{i=1}^\infty \PP[A_i] = \infty$. Then
\begin{align} 
 \limsup_{n \to \infty} \frac{(\sum_{k=1}^n \PP[A_k])^2}{\sum_{i,k=1}^n \PP[A_i A_k]}
 & = \limsup_{n \to \infty} \frac{\sum_{1 \le i < k \le n} \PP[A_i]\PP[A_k]}{\sum_{1 \le i < k \le n}\PP[A_iA_k]} \label{yan-eq}
\end{align}
\end{thm}
\begin{proof} Define sequences $s_n, t_n, b_n$ and $c_n$ as follows:
\begin{align*}
s_n &= \sum_{k=1}^n \PP[A_k] &
    t_n &= \sum_{1 \le i < k \le n} \PP[A_i]\PP[A_k] \\
b_n &= \sum_{i,k=1}^n \PP[A_i A_k] &
    c_n  &= \sum_{1 \le i < k \le n}\PP[A_iA_k]
\end{align*}
Since $\sum_{i=1}^\infty \PP[A_i] = \infty$, we have $s_n = o(s_n^2)$.
Hence, as $2t_n \le s_n^2 \le 2t_n + s_n$, $\lim_{n\to \infty} s_n^2/2t_n = 1$.
By inequality (\ref{kochen-stone-ineq}), $s_n^2 \le (1 + o(1)) b_n$.
Hence, as $2 c_n \le b_n = 2 c_n + s_n$, $\lim_{n\to\infty} b_n/2 c_n = 1$.
It follows that $\limsup_{n\to\infty} s_n^2/b_n = \limsup_{n \to \infty} t_n/c_n$, which is what we wish to prove.
\end{proof}

Let $(q_n)_{n=1}^\infty$ be any non-decreasing sequence of elements of the open unit interval $(0, 1)$, let $q = \lim_{n \to \infty} q_n$ and let $A_i$ be the event that a uniformly random member of the unit interval $[0, 1]$ lies in $[0, q_i]$.
Then $\PP[A_i] = q_i$ and $\PP[A_i~\mbox{i.o.}] = q$.
Moreover, $\PP[A_i A_k] = \PP[A_i]$ for $i < k$.
Let us define $u_n, v_n$ and $w_n$ as follows:
\begin{align}
    u_n &= q_1\sum_{k= 2}^n q_k + q_2\sum_{k=3}^n q_k + \ldots + q_{n-1}q_n\\
    v_n &= (n-1)q_1 + (n-2)q_2 + \ldots + q_{n-1}\\
    w_n &= \frac{u_n}{v_n}
\end{align}
Then, rearranging the terms in the sums on the right-hand side of equation (\ref{yan-eq}), we find that, given equation (\ref{yan-eq}), the inequality (\ref{kochen-stone-ineq}) is equivalent to
$q \ge \limsup_{n \to \infty} w_n$. The following lemma implies that equality holds in the Kochen-Stone inequality for any sequence of events $(A_i)_{i=1}^\infty$ constructed in this way.


\begin{lem}\label{lma:lim-wn} Let $w_n$ and $q$ be as above. Then $w_n \to q$ as $n \to \infty$.
\end{lem}
\begin{proof}
We have:
\begin{align}
    \sum_{k=i}^nq_k &= q\left(n - i + 1 - \frac{1}{q}\sum_{k=i}^n (q - q_k) \right)
\end{align}
Let us write $\sigma_i^j$ for $\sum_{k=i}^j(q - q_k)$.
From the above, multiplying by $q_{i-1}$ and summing for $i$ from $2$ to $n$, we have:
\begin{equation}
    u_n = qv_n - \sum_{i=2}^{n}  q_{i-1}\sigma_i^n
\end{equation}
Define the sequence $r_n$ by:
\begin{align}
    r_n &= \frac{\sum_{i=2}^{n}  q_{i-1}\sigma_i^n}{v_n} 
\end{align}
We claim that $r_n \to 0$ as $n \to \infty$ so that $w_n = u_n/v_n = (q - r_n) \to q$, which is what we have to prove.
So given $\varepsilon > 0$, let $\varepsilon_0 = \frac{q_1}{4}\varepsilon$ and choose $N$ such that for all $n > N$,we have $q-q_n < \varepsilon_0$.
Define $C$ by:
\begin{equation}
    C = q_1\sigma_2^N + q_2 \sigma_3^N + \ldots q_{N-1}\sigma_N^N
\end{equation}
Then, for $n > N$:
\begin{align}
    r_n &= \frac{C}{v_n} +
    \frac{q_1\sigma_{N+1}^n + \ldots + q_N\sigma_{N+1}^n + q_{N+1}\sigma_{N+2}^n + \ldots + q_{n-1}\sigma_n^n}{(n-1)q_1 + (n-2)q_2 + \ldots + q_{n-1}} \label{eqn-req-c-over-v-plus}\\
    &\le \frac{C}{v_n} + \frac{(n-1)\sigma_{N+1}^n}{q_1((n- 1) + (n - 2) + \ldots + 1)} \label{eqn-r-bound-a}\\
    &\le \frac{C}{v_n} + \frac{1}{q_1}\cdot\frac{(n-1)(n-N)\varepsilon_0}{\frac{1}{2}(n-1)n} \label{eqn-r-bound-b}\\
    &\to \frac{2}{q_1}\varepsilon_0 = \frac{\varepsilon}{2}\mbox{ as $n \to \infty$}
\end{align}
where in equation (\ref{eqn-req-c-over-v-plus}) we have expanded the denominator of the second fraction using the definition of $v_n$and where the bounds (\ref{eqn-r-bound-a}) and (\ref{eqn-r-bound-b}) are obtained using the facts that $0 < q_1 \le q_i < 1$ and $\sigma_{N+m}^n \le \sigma_{N+1}^n \le (n - N) \varepsilon_0$.
Hence we can choose $M > N$, such that for $n > M$, we have $|r_n - \frac{\varepsilon}{2}| < \frac{\varepsilon}{2}$, giving $r_n < \varepsilon$. Hence $r_n \to 0$ as $n \to \infty$.
\end{proof}



Recall that a \emph{Specker sequence} (see \cite{Specker49} or \cite{Troelstra-van-Dalen88}) is a computable, monotone increasing, bounded sequence of rationals whose limit is not a computable real number\footnote{A real number is said to be \emph{computable} if it is the limit of a computable sequence of rationals, with a computable rate of convergence (see \cite{Weihrauch(00)}).}. 

\begin{thm}
There is a sequence of events $(A_i)_{i=1}^\infty$
and a computable function $\omega : \NN^+ \to \NN^+$ such that, for any $N$,
\[ \left(\sum_{i=1}^{\omega(N)} \PP[A_i]\right) \geq N \]
but for which there is no computable function $\phi : \NN \times \NN^+ \to \NN^+$ satisfying
\begin{equation}
\forall m, \ell\exists n \in [m, \phi(m, \ell)] \left( \PP\left[\bigcup_{i=m+1}^n A_i\right] + \frac{1}{2^\ell} \geq \limsup_{j \to \infty} \frac{(\sum_{k=1}^j \PP[A_k])^2}{\sum_{i,k=1}^j \PP[A_i A_k]} \right)
\label{eqn:bound}
\end{equation}
\end{thm}
\begin{proof} Take $A_i$ to be the event that a uniformly random element of $[0, 1]$ lies in the interval $[0, q_i]$ where the $q_i \in (0, 1)$ form a Specker sequence with limit $q$ (note that $q < 1$, since $q$ is not a computable real). 
Also $\sum_{i=1}^\infty q_i$ diverges with the rate of divergence given by the computable function $\omega(N) = \lceil\frac{N}{q_1}\rceil$.
Then, for $n > m$, we have:
\begin{equation}
    \PP\left[\bigcup_{i=m+1}^n A_i\right] = \PP[A_n] = q_n
\end{equation}
and, by Lemma~\ref{lma:lim-wn} and the discussion preceding it, we have:
\begin{equation}
    \limsup_{j \to \infty} \frac{(\sum_{k=1}^j \PP[A_k])^2}{\sum_{i,k=1}^j \PP[A_i A_k]} = \lim_{j \to \infty} \frac{(\sum_{k=1}^j \PP[A_k])^2}{\sum_{i,k=1}^j \PP[A_i A_k]} = q.
\end{equation}
 Assume $\phi$ satisfies (\ref{eqn:bound}), i.e.
\begin{equation}
\forall m, \ell \exists n \in (m, \phi(m, \ell)] \left( q_n + \frac{1}{2^\ell} \geq q \right)
\end{equation}
Hence
\begin{equation}
\forall \ell \exists n \leq \phi(0, \ell) \left( |q - q_n| \leq \frac{1}{2^\ell} \right)
\end{equation}
implying that, $g(\ell) = \phi(0, \ell)$ is a rate of convergence for the Specker sequence $q_n$. But since $q = \lim_{n\to\infty} q_n$ is not a computable real, it cannot be approximated by a sequence with a computable rate of convergence. It follows that $g$, and hence also $\phi$, is not a computable function.
\end{proof}

\section{Optimality of the Estimates}\label{sec-optimal}

It is easy to argue that the numeric bounds in Theorem \ref{thm:first-bcl-q} are best possible. Indeed, let $(A_i)_{i=1}^\infty$ be a sequence of mutually exclusive events. In this case, the first inequality in (\ref{first-lemma-proof}) is actually an equality and hence the given estimate is optimal.

The estimate given by Theorem~\ref{thm:second-bcl-q} is also the best possible amongst estimates that do not depend on $n$. To see this, consider the probability space whose outcomes are functions $\alpha \colon \NN^+ \to \{1,\ldots, k\}$ representing an infinite sequence of throws of a fair $k$-sided die. Let $A_i$ be the event $\alpha(i) = k$, so that $\PP[A_i] = 1/k$. Clearly $\sum_{i=1}^\infty \PP[A_i]$ diverges with (optimal) rate $\omega(N) = k N$. We have that
\begin{equation}
    \PP\left[ \bigcap_{i=n}^{k (n + N - 1)} \overline{A}_i \right] = \left( 1 - \frac{1}{k} \right)^{k (n + N - 1) - n + 1} 
\end{equation}
This is a decreasing function of $n$, so the worst case for our estimate is when $n = 1$, but in that case we have:
\begin{equation}
   \PP\left[ \bigcap_{i=n}^{k N} \overline{A}_i \right] =  \left( 1 - \frac{1}{k} \right)^{k N} \to e^{-N} \quad \mbox{as $k \to \infty$}
\end{equation}
showing that the bounds in (\ref{2nd-bc-proof-eq1}) and (\ref{2nd-bc-proof-eq2}) are tight.

Unlike the quantitative proofs of the First and Second Borel-Cantelli lemmas, where we have argued above that the numerical bounds obtained have optimal rate of growth, in the case of the quantitative Erd\H{o}s-R\'{e}nyi theorem (Theorem \ref{erdos-renyi-thm}), it is plausible that the bounds we have achieved are optimal for the given formulation of assumption (\ref{erdos-assumption2}) (see remark after Theorem \ref{erdos-renyi-thm}), but we have been unable to prove that. 

In the case of the quantitative Kochen-Stone theorem (Theorem \ref{quant-kochen-stone-thm}), we have had to reformulate the original theorem using Tao's notion of metastability. This means that the numerical bound given also depends on a new functional input $g \colon \NN \to \NN$. This makes an asymptotic analysis of the growth rate of the bound much more difficult. It is again unclear whether the bound we have obtained has optimal dependency on any of the parameters $\omega, g, m$ and $\ell$. We leave this for future work. It is worth emphasising, however, that all our numerical bounds are uniform in (i.e. independent of) the actual sequence of events $(A_i)_{i=1}^\infty$.


\section*{Acknowledgements}

We would like to thank Iosif Pinelis for the proof of Theorem~\ref{thm:yan}, and the anonymous referees for their constructive feedback, corrections and suggestions.


%
%
%
\bibliographystyle{jloganal}
%


%

\end{document}